\documentclass[11pt]{article}
\usepackage{amsfonts,amsmath,latexsym,color,epsfig}
\date{}
    \newcommand{\C}{\mathcal{C}}
    \newcommand{\D}{\mathcal{D}}
    \newcommand{\M}{\mathcal{M}}
    \newcommand{\F}{\mathcal{F}}
    \newcommand{\G}{\mathcal{G}}

    \newcommand{\eps}{\varepsilon}
    \newcommand{\pr}{{\rm {\bf Pr}}}
    \newcommand{\VC}{{\mbox{VC-dim}}}

\newtheorem{theorem}{Theorem}
\newtheorem{lemma}[theorem]{Lemma}
\newtheorem{corollary}[theorem]{Corollary}
\newtheorem{claim}[theorem]{Claim}
\newtheorem{definition}[theorem]{Definition}

\newtheorem{problem}[theorem]{Problem}

\def\qed{\ifhmode\unskip\nobreak\fi\quad\ifmmode\Box\else$\Box$\fi}

\title{Shattered matchings in intersecting hypergraphs}

\author{
{\sl Peter Frankl}\thanks{R\'enyi Institute, P.O.Box 127 Budapest, 1364 Hungary; \texttt{peter.frankl@gmail.com}.} \thanks{MIPT, Moscow, partially supported by the Ministry of Education and Science of the Russian Federation in the framework of MegaGrant no 075-15-2019-1926.}
\and
{\sl J\'{a}nos Pach}\thanks{R\'enyi Institute, P.O.Box 127 Budapest, 1364 Hungary; \texttt{pach@cims.nyu.edu}, partially supported by NKFIH \'Elvonal (Frontier) program KKP 133864.} \thanks{IST Austria, Vienna, partially supported by Austrian Science Fund (FWF), grant Z 342-N31.}
\thanks{MIPT, Moscow, partially supported by the Ministry of Education and Science of the Russian Federation in the framework of MegaGrant no 075-15-2019-1926.}}

\date{}

\begin{document}

\maketitle

\begin{abstract}
Let $X$ be an $n$-element set, where $n$ is even. We refute a conjecture of J. Gordon and Y. Teplitskaya, according to which, for every maximal intersecting family $\mathcal{F}$ of $\frac{n}2$-element subsets of $X$, one can partition $X$ into $\frac{n}2$ disjoint pairs in such a way that no matter how we pick one element from each of the first $\frac{n}2 - 1$ pairs, the set formed by them can always be completed to a member of $\mathcal{F}$ by adding an element of the last pair.

The above problem is related to classical questions in extremal set theory. For any $t\ge 2$, we call a family of sets $\F\subset 2^X$ {\em $t$-separable} if for any ordered pair of elements $(x,y)$ of $X$, there exists $F\in\F$ such that $F\cap\{x,y\}=\{x\}$.   For a fixed $t, 2\le t\le 5$ and $n\rightarrow\infty$, we establish asymptotically tight estimates for the smallest integer $s=s(n,t)$ such that every family $\F$ with $|\F|\ge s$ is $t$-separable.
\end{abstract}

\section{Introduction}\label{intro}

Given an $n$-element set $X$, a family $\F\subset 2^X$ is called {\em intersecting} if any two members of $\F$ have nonempty intersection. In their seminal work~\cite{ErKR61}, P. Erd\H os, C. Ko, and R. Rado determined the maximum size of an intersecting family $\F$ of $k$-element subsets of $X$, for all $k\le \frac{n}2$. In particular, if $n$ is even and $k=\frac{n}2$, they proved that $|\F|\le\frac12\binom{n}{\frac{n}2}$, where equality holds for every {\em maximal} (that is, non-extendable) intersecting family of $\frac{n}2$-element subsets of $X$.
\smallskip

Motivated by a problem from mathematical finance, J. Gordon and Y. Teplitskaya~\cite{GoT19} made the following conjectures:
\smallskip

\noindent {\bf Conjecture A.} If $n$ is {\em even}, then for any {\em maximal intersecting} family $\F$ of $\frac{n}2$-element subsets of $X$, there exists a {\em perfect matching} $\{x_1,x_2\},\{x_3,x_4\},\ldots,\{x_{n-1},x_n\}$ with the property that no matter how we select one element from each of the first $\frac{n}{2}-1$ pairs, together with $x_{n-1}$ or $x_n$, they always form a member of $\F$.

\noindent {\bf Conjecture B.} If $n$ is {\em odd}, then for any {\em maximal intersecting} family $\F$ consisting of $\frac{n-1}2$-element and $\frac{n+1}2$-element subsets of $X$, there exists a matching
$\{x_1,x_2\},\{x_3,x_4\},\ldots,\{x_{n-2},x_{n-1}\}$ with the property that no matter how we select one element from each pair, together with the last element $x_n\in X$, they always form a member of $\F$.
\smallskip

Gordon and Teplitskaya verified these conjectures for $n\le 6$.

In this note, we disprove the above conjectures for all $n\ge 14$. For $n$ even, we will establish a more general result which contradicts Conjecture A in a strong way. For odd $n$, the problem will be settled using the even case. To formulate our first result, we need to agree on some terminology.
\smallskip

An unordered collection $\{x_1,x_2\},\{x_3,x_4\},\ldots,\{x_{2k-1}x_{2k}\}$ of pairwise disjoint $2$-element subsets of $X$ is called a {\em matching} of size $k$\; ($2k\le n=|X|$). If we pick one element from each pair, the $k$-element set formed by them is called a {\em snake} with respect to this matching. A family of subsets of $X$ is said to be {\em intersecting} if any two of its members have nonempty intersection.

\begin{definition}\label{shattered}
A {\em matching} $\{x_1,x_2\},\{x_3,x_4\},\ldots,\{x_{2k-1}x_{2k}\}\subset X$ is said to be {\em shattered} by a family $\F\subset 2^X$ if for every snake $S$ with respect to this matching, there exists $F\in\F$ such that $F\cap\{x_1,x_2,\ldots,x_{2k}\}=S$.
\end{definition}

Our main result is the following.

\begin{theorem}\label{thm:main}
Let $X$ be an $n$-element set, where $n\ge 28$ is even. Let $k(n)$ denote the largest integer $k$ such that for every maximal intersecting family of $\frac{n}2$-element subsets of $X$, there exists a shattered matching of size $k$. Then we have
$$\frac{n}{4}\le k(n) \le \frac{n}{2}-\frac12\log_2n + 1.$$
\end{theorem}

The fact that the upper bound is smaller than $\frac{n}{2}-1$  shows that the Gordon-Teplitskaya conjecture is not true if $n$ is large enough.
\smallskip

The above question is closely related to a classical result from extremal set theory~\cite{Fr95, FrT18, GeP19}.  A family $\F$ of subsets of $X$ is said to {\em shatter a set} $A$ if for every subset $B\subset A$, there is $F\in \F$ with $F\cap A=B$. It was shown by Vapnik and Chervonenkis and a little later, independently, by Sauer and Shelah that if $\F$ is large, then there is a large subset $A\subset X$ shattered by $\F$. The size of the largest shattered subset of $X$ is called the {\em Vapnik-Chervonenkis dimension} of $\F$, and is denoted by $\VC(\F)$. This notion plays a central role in statistics, learning theory, discrete and computational geometry, and elsewhere.

More precisely, the following is true.

\begin{theorem}\label{sauer}
{\rm (Sauer~\cite{Sa72}, Shelah~\cite{Sh72}, Vapnik-Chervonenkis~\cite{VaCh71})}  Let $|X|=n$ and let $\F$ be a family of subsets of $X$ with $|\F|>\sum_{i=0}^{k-1}{n\choose i}$.

Then there is a $k$-element set $A\subset X$ shattered by $\F$, i.e., $\VC(\F)\ge k$. This bound is tight.
\end{theorem}

Obviously, if $A$ is shattered and $|A|$ is even, then any perfect matching of $A$ (that is, any partition of $A$ into $2$-element subsets) is a shattered matching. Using this idea, one can easily obtain the lower bound $k(n)\ge (\frac14-o(1))n$, which is only slightly weaker than the bound in Theorem~\ref{thm:main}. However, this proof only uses that $\F$ has many members, without taking into account the assumption that $\F$ is a maximal intersecting family.

In the spirit of Theorem~\ref{sauer}, we can ask how large $\F$ must be in order to guarantee the existence of a shattered matching of size $k$.

\begin{problem}\label{problem2}
Let $n, k$ be positive integers, $n\ge 2k$. Determine or estimate the smallest number $p=p(n,k)$ such that for every family $\F$ of at least $p$ subsets of an $n$-element set, there is a matching of size $k$ shattered by $\F$.
\end{problem}

Obviously, if $\F$ shatters a $2k$-element set $A$, then any partition of $A$ into $2$-element sets is a shattered matching of size $k$. Thus, Theorem~\ref{sauer} immediately implies that
$$p(n,k)\le 1+ \sum_{i=0}^{2k-1}{n\choose i}\;\;\;\;\; \mbox{ for every }\; k\ge 1.$$
This bound is tight for $k=1$, and we will see that its order of magnitude is best possible for any $k$, as $n\rightarrow\infty$. See Corollary~\ref{cor1}.

\begin{definition}\label{separable}
A family $\F\subset 2^X$ is said to be {\em $t$-separable} if there is a $t$-element subset $T\subset X$ such that for every ordered pair $x,y\in T,\; x\neq y$, there exists $F\in\F$
such that $F\cap\{x,y\}=x$.

For any $n\ge t\ge 2,$ let $s(n,t)$ denote the smallest number $s$ with the property that every family $\F$ of at least $s$ subsets of an $n$-element set is $t$-separable.
\end{definition}

If $\{x_1,x_2\},\ldots,\{x_{2k-1}x_{2k}\}$ is a matching of size $k$ shattered by $\F$, then $\F$ is $2k$-separable, as the set $T=\{x_1,\ldots,x_{2k}\}$ satisfies the above requirements. Therefore, we have
$$p(n,k)\ge s(n,2k).$$

The problem of determining or estimating $s(n,t)$ appears to be a nontrivial task of independent interest. We prove the following.

\begin{theorem}\label{sep}
Let $n\ge t\ge 2$, and let $X$ be an $n$-element set. Let $s(n,t)$ denote the smallest number $s$ with the property that every family $\F\subset 2^X$ with $|\F|\ge s$ is $t$-separable.

{\rm (i)}\;\;\; For $t=2$, we have $s(n,2)=n+2$.

{\rm (ii)}\;\;  For $t=3$, we have $s(n,3)=\lfloor\frac{n^2}{4}\rfloor+n+2$.

{\rm (iii)}\;   For $t=4$ or $5$, we have $s(n,t)=(\frac{n}{t-1})^{t-1} + \Theta(n^{t-2})$.

{\rm (iv)}\;\;  For $t\ge 6$, we have  $(\frac{n}{t-1})^{t-1}<s(n,t)\le 1+\sum_{i=0}^{t-1}{n\choose i}.$
\end{theorem}

A more precise form of parts (iii) and (iv) is stated and proved as Theorem 6' in Section~\ref{section3}.
\smallskip

This note is organized as follows. In the next section, we prove Theorem~\ref{thm:main}. In Section~\ref{section3}, we study the function $s(n,t)$, and we establish Theorem~\ref{sep} (and Theorem~6'). The last section contains some open problems and related results.

\section{Shattered matchings---Proof of Theorem~\ref{thm:main}}\label{section2}

Throughout this note, let $A\sqcup B$ denote the {\em disjoint union} of the sets $A$ and $B$. For any $F, Y\subset X$, we call $F\cap Y$ the {\em trace of $F$ on $Y$.} For a family $\F\subset 2^X$ and $Y\subset X$, the set of traces $F\cap Y$ over all $F\in\F$ is denoted by $\F_{|Y}$.

\medskip

\noindent{\bf Proof of the lower bound.} Let $\F\subset 2^X$ be a maximal intersecting family of $\frac{n}2$-element subsets of $X$. By the maximality of $\F$, if an $\frac{n}2$-element set $Y\subset X$ does not belong to $\F$, then its complement $\overline{Y}=X\setminus Y$ does.

Let $M=\{x_1,x_2\}\sqcup\ldots\sqcup\{x_{2k-1},x_{2k}\}$ be any matching of size $k$ in $X$. Extend it to a perfect matching of $X$ by adding a perfect matching of the remaining $n-2k$ elements: $M'= \{x_{2k+1},x_{2k+2}\}\sqcup\ldots\sqcup\{x_{n-1},x_{n}\}$. If $M$ is not a shattered matching, then we can choose a snake $S$ of $M$ such that no member $F\in\F$ intersects $\{x_1,\ldots,x_{2k}\}$ precisely in the elements of $S$. Suppose without loss of generality that $\{x_1, x_3,\ldots, x_{2k-1}\}$ is such a snake. This implies that
$$ \{x_1, x_3,\ldots, x_{2k-1}\}\sqcup\{x_{2k+2-\eps(1)}, x_{2k+4-\eps(2)},\ldots,x_{n-\eps(\frac{n}2-k)}\}$$
does not belong to $\F$ for any $\eps(i)\in\{0,1\},\; 1\le i\le \frac{n}2-k$. Hence, the complement of this set,
$$ \{x_2, x_4,\ldots, x_{2k}\}\sqcup\{x_{2k+1+\eps(1)}, x_{2k+3+\eps(2)},\ldots,x_{n-1+\eps(\frac{n}2-k)}\}\in\F,$$
for every $\eps(i)\in\{0,1\}$. This means, by definition, that $M'$ is a shattered matching of size is $\frac{n}2-k$.

Thus, either there is a shattered matching $M$ of size $k$, or a shattered matching $M'$ of size $\frac{n}2-k$. \hfill  $\Box$
\medskip

\noindent{\bf Proof of the upper bound.}  Two matchings of the same size in $X$ are considered identical if they differ only in the order of pairs. We need some simple facts.

\begin{claim}\label{claim1}
The number of matchings of size $k$ in $X$ is smaller than $\frac{n!}{k!2^k}$.
\end{claim}

Indeed, with each permutation $(x_1,x_2,\ldots,x_n)$ associate the matching $\{x_1,x_2\},\ldots,\{x_{2k-1},x_{2k}\}$, and note that we obtain every matching at least $k!2^k$ times.
\smallskip

Now we {\em randomly} generate a maximal intersecting family $\F\subset 2^X$ consisting of $\frac{n}2$-element subsets of $X$, as follows. From each of the $\frac12\binom{n}{n/2}$ unordered pairs $(Y,X\setminus Y)$ with $|Y|=n/2$, we select either $Y$ or $X\setminus Y$, independently with probability $\frac12$. Let $\F$ consist of all the selected sets. Obviously, any two members of $\F$ have nonempty intersection and $\F$ is maximal with respect to this property.

Fix a matching (partition) $M=\{x_1,x_2\}\sqcup\ldots\sqcup\{x_{2k-1},x_{2k}\}$ of size $k$ in $X$. Let $S$ be a snake with respect to $M$. We say that $S$ is {\em carved out of $M$} by a family $\F\subset 2^X$ if there exists $F\in \F$ whose trace on $\{x_1,\ldots,x_{2k}\}$ is $S$, i.e., if we have $F\cap\{x_1,\ldots,x_{2k}\}=S$. With a slight abuse of notation, we write $M\setminus S$ for the set $\{x_1,\ldots,x_{2k}\}\setminus S$
which is also a snake with respect to the matching $M$.

In the sequel, for convenience, we write $2\ell$ for $n-2k$, so that $k+\ell=\frac{n}2$.

\begin{claim}\label{claim2}
Let $M$ be a fixed matching of size $k=\frac{n}2-\ell$ in $X$. For any snake $S$ with respect to $M$, the probability that $S$ is {\em not carved out} of $M$ by the randomly generated family $\F$ is equal to $2^{-\binom{2\ell}{\ell}}$.
\end{claim}

To see this, it is enough to notice that if $Y\cap\{x_1,\ldots,x_{2k}\}$ is not equal to $S$, nor to $M\setminus S$, then it does not matter which set we select from the pair $(Y,X\setminus Y)$, it can not separate $S$ from $M$. Therefore, it is enough to consider the $\binom{2\ell}{\ell}$ pairs of complementary sets $(Y,X\setminus Y),$ where $Y=S\sqcup T$ for some $\ell$-element subset $T\subset X\setminus M$. For each of these pairs, we have to select the set $X\setminus Y$ to be contained in $\F$, otherwise $S$ will be carved out by $Y$ and, hence, by $\F$. This proves Claim~\ref{claim2}. The probability $2^{-\binom{2\ell}{\ell}}$ may appear to be tiny, but for a fixed $\ell$ it is bounded away from $0$, as $n\rightarrow\infty$.
\smallskip

Notice that if $S$ and $S'$ are two distinct snakes with respect to $M$ and $S'\not=M\setminus S$, then the events that ``$S$ {\em is carved out of} $M$'' and ``$S'$ {\em is carved out of} $M$'' are independent, because they depend on completely different random choices. On the other hand, for $S'=M\setminus S$, we have
\begin{align*}
&\pr[\mbox{at least one of } S \mbox{ and } M\setminus S \mbox{ is not carved out of } M]\\
&=\pr[\mbox{precisely one of } S \mbox{ and } M\setminus S \mbox{ is not carved out of } M]\\
&=2\cdot\pr[S \mbox{ is not carved out of } M]= 2\cdot2^{-\binom{2\ell}{\ell}}.
\end{align*}

The number of unordered pairs of snakes $(S,M\setminus S)$ is $2^{k-1}$. Hence,
\begin{align*}
&\pr[M \mbox{ is shattered}]\\
&=\pr[\mbox{every pair of snakes } (S,M\setminus S) \mbox{ are carved out of } M]\\
&=(1-2\cdot2^{-\binom{2\ell}{\ell}})^{2^{k-1}}<\exp(-2^{k-\binom{2\ell}{\ell}}).
\end{align*}

Combining this with Claim~\ref{claim1}, we obtain that
\begin{align*}
&\pr[\mbox{there exists a shattered matching of size } k]\\
&\leq\sum_{M}\pr[M \mbox{ is shattered}]\\
&<\frac{n!}{k!2^k}\exp(-2^{k-\binom{2\ell}{\ell}})<\exp(n\ln n-2^{\frac{n}{2}-\ell-\binom{2\ell}{\ell}}).
\end{align*}
\noindent To conclude, it is enough to show that the right-hand side of this inequality is smaller than 1, that is, $2^{\frac{n}{2}-\ell-\binom{2\ell}{\ell}}>n\ln n$ holds, provided that $k\ge \frac{n}2-\frac12\log_2n +1$.

According to the last condition, $\ell+1\le\frac12\log_2n$, which implies that
$$2^{\frac{n}{2}-\ell-\binom{2\ell}{\ell}}\ge 2^{\frac{n}{2}-4^\ell}\ge 2^{\frac{n}{4}}>n\ln n,$$
if $n\ge28$.
This completes the proof of the upper bound and, hence, Theorem~\ref{thm:main}. \hfill $\Box$
\smallskip

It is easy to verify using the above estimates that the probability that there exists a shattered matching of size $k=\frac{n}2-1$ is smaller than $1$, for every $n\ge 14$. Therefore, in these cases, Conjecture A of Gordon and Teplitskaya fails.
\smallskip

Next, we turn to Conjecture B.

\begin{corollary}\label{odd}
Let $X$ be an $n$-element set, where $n\ge 15$ is odd.

There is a maximal intersecting family $\F\subset 2^X$ such that $|F|=\frac{n-1}2$\; or $\frac{n+1}2$ for every $F\in\F$, and the following condition is satisfied. There do not exist $y\in X$ and a perfect matching $M$ of $X\setminus\{y\}$ with the property that adding $y$ to  every snake with respect to $M$, we get a member of $\F$.
\end{corollary}

\noindent{\bf Proof.} Let $|X|=n=2k+1$, let $\F\subset 2^X$ be a maximal intersecting family, and suppose that every member of $\F$ has $k$ or $k+1$ elements. By the maximality of $\F$, for each $k$-element subset $Y\subset X$, either $Y$ or $X\setminus Y$ belongs to $\F$. Therefore, we have $|\F|=\binom{n}{k}$.
\smallskip

Fix a $2k$-element subset $V\subset X$, and denote the unique element of $X\setminus V$ by $x$.
According to the remark after the proof of Theorem~\ref{thm:main}, we can choose a maximal intersecting family $\G$ of $\frac12\binom{2k}{k}$ $k$-element subsets of $V$ such that $\G$ does not shatter any matching of size $k-1$ in $V$. Let
$$\F=\G\bigcup\{\,Y\cup\{x\}\,|\,Y\in\G\,\}\bigcup\{\,U\,|\, U\subset V, |U|=k+1\, \}.$$
Obviously, $\F$ is a {\em maximal intersecting} family consisting of $k$-element and $(k+1)$-element subsets of $X$.

We claim that $\F$ meets the requirements of Corollary~\ref{odd}.
Suppose for contradiction that there are $y\in X$ and a perfect matching $M$ of $X\setminus\{y\}$ such that every snake with respect to $M$ can be extended to a member of $\F$ by adding $y$. We distinguish two cases.
\smallskip

Suppose first that $y=x$. Then $X\setminus\{y\}=X\setminus\{x\}=V$, and $M$ is a partition of $V$ into $2$-element sets:
$$\{x_1,x_2\}\sqcup\{x_3,x_4\}\sqcup\ldots\sqcup\{x_{2k-1},x_{2k}\}.$$
It follows from the definition of $\G$ that only one of the snakes with respect to $M$, $\{x_1,x_3,\ldots,x_{2k-1}\}$ or $\{x_2,x_4,\ldots,x_{2k}\}$ belongs to $\G$. Hence, only one of the ``extended'' snakes $\{x_1,x_3,\ldots,x_{2k-1},y\}$ or $\{x_2,x_4,\ldots,x_{2k},y\}$ belongs to $\F$. The other one does not, contradicting our assumption.
\smallskip

Suppose next that $y\neq x$, and let $M$ be the partition (perfect matching) of $X\setminus\{y\}$ with the above property,
$$X\setminus\{y\}=\{x_1,x_2\}\sqcup\{x_3,x_4\}\sqcup\ldots\sqcup\{x_{2k-1},x_{2k}\}.$$
We can assume, by symmetry, that $x=x_{2k}$. It follows from the definition of $\G$ that the ``partial'' matching $M'$ of $M$, $$\{x_1,x_2\}\sqcup\{x_3,x_4\}\sqcup\ldots\sqcup\{x_{2k-3},x_{2k-2}\},$$
is not shattered by $\G$.

We can assume without loss of generality that the snake $\{x_1,x_3,\ldots,x_{2k-3}\}$ with respect to $M'$ is {\em not} the trace of any member of $\G$ on $V\setminus\{x_{2k-1},y\}$.
In particular, we have $\{x_1,x_3,\ldots,x_{2k-3},y\}\not\in\G$. By the definition of $\F$, this implies that
$$\{x_1,x_3,\ldots,x_{2k-3},y,x\}\not\in\F.$$
However, this means that the snake $\{x_1,x_3,\ldots,x_{2k-3},x\}$ with respect to the matching $M$, cannot be extended to a member of $\F$ by adding $y$, contradiction. This completes the proof of the corollary. \hfill $\Box$

\section{Separable families---Proof of Theorem~\ref{sep}}\label{section3}

We start with a construction of non-$t$-separable families. Let $x_1,\ldots,x_n$ be the elements of a set $X$, listed in an arbitrary order. The set-system $$\C=\{\emptyset,\{x_1\},\{x_1,x_2\},\ldots\{x_1,\ldots,x_n\}\}$$
is called a {\em maximal chain} on $X$.

\begin{lemma}\label{constr}
For $n\ge t\ge 2$, consider a partition of an $n$-element set $X$ into $t-1$ parts, $X=X_1\sqcup \ldots\sqcup X_{t-1}$. Fix a maximal chain $\C_i$ on each $X_i$.

Then the family
$$\D(X_1,\ldots,X_{t-1})=\{C_1\sqcup\ldots\sqcup C_{t-1}\; :\: C_i\in\C_i \mbox{ for } i=1,\ldots,t-1\}$$
is not $t$-separable.
\end{lemma}

\noindent{\bf Proof.} Suppose for contradiction that $\D(X_1,\ldots,X_{t-1})$ is $t$-separable, that is, there exists a $t$-element subset $T\subset X$ satisfying the conditions in Definition \ref{sep}. By the pigeonhole principle, there is an $X_i,\; 1\le i\le t-1$ which contains at least two elements, $x,y\in T$. Suppose without loss of generality that $x$ precedes $y$ in the order that defines the chain $\C_i$ on $X_i$. Then $\C_i$ and, hence, $\D(X_1,\ldots,X_{t-1})$ has no member which contains $y$, but not $x$.                            \hfill $\Box$
\smallskip

Similar constructions involving direct products of {\em chains} can be found, e.g., in~\cite{An13, AnS05}.

\begin{corollary}\label{cor1}
Let $n\ge t\ge 2$, and let $X_1\sqcup \ldots\sqcup X_{t-1}$ be a partition of an $n$-element set into $t-1$ parts, as equal as possible. That is, we have $|X_i|=\lfloor\frac{n}{t-1}\rfloor$ or $\lceil\frac{n}{t-1}\rceil$ for $1\le i\le t-1$. Then the smallest number $s=s(n,t)$ with the property that every family of at least $s$ subsets of $X$ is $t$-separable, satisfies
$$s(n,t)>\prod_{i=1}^{t-1}(|X_i|+1)>(\frac{n}{t-1})^{t-1}.$$
Consequently, for every $k\le\frac{n}{2}$, we have
$$p(n,k)\ge s(n,2k)>(\frac{n}{2k-1})^{2k-1}.$$  \hfill $\Box$
\end{corollary}

For $t=2$, the first part of Corollary~\ref{cor1} implies that $s(n,2)\ge n+2$. On the other hand, it follows from Theorem~\ref{sauer} that if $|\F|>\binom{n}1+\binom{n}0=n+1$ for a family $\F$ of subsets of an $n$-element set $X$, then $\F$ shatters a $2$-element subset of $X$, hence, $\F$ is $2$-separable. Thus, $s(n,2)=n+2$, which proves part (i) of Theorem~\ref{sep}.
\smallskip

In what follows, we use the ``arrow'' notation proposed by Hajnal (see~\cite{Bo72, ErHMR84}). We write $$(n,m)\rightarrow (a,b)$$ if the following statement is true: For any family of $m$ subsets $\F=\{F_1,\ldots,F_m\}$ of an $n$-element set $X$, there is $T\subset X$ with $|T|=a$ such that the family of {\em traces} $\F_{|T}=\{F_1\cap T,\ldots,F_m\cap T\}$ has at least $b$ distinct members.

Using this notation, Theorem~\ref{sauer} can be reformulated as
$$(n, 1+\sum_{i=0}^{k-1}\binom{n}{i})\rightarrow (k, 2^k).$$

A family $\F$ is called {\em downward closed} if for any $F\in\F$ and $G\subset F$, we have $G\in \F$. A crucial property of the arrow relation was established by the first author~\cite{Fr83}.

\begin{lemma}\label{downward}
{\rm (Frankl~\cite{Fr83})}  The relation $(n,m)\rightarrow (a,b)$ holds if and only if for any {\em downward closed} family $\F$ of $m$ subsets of a set $X$ with $|X|=n$, there is $T\subset X$ with $|T|=a$ such that $|\F_{|T}|\ge b$.
\end{lemma}

In other words, in order to show that there is $T\subset X$ with $|T|=a$ and $|\F_{|T}|\ge b$, it is sufficient to verify it for downward closed families.

\begin{lemma}\label{lem2}
Suppose that $(n,m)\rightarrow (t,2^t-2^{t-2}+1)$ holds.

Then every family $\F$ of subsets of an $n$-element set with $|\F|\ge m$ is $t$-separable.
\end{lemma}

\noindent{\bf Proof.}  Assume that $\F\subset 2^X$ satisfies the above condition, and let $T$ be a $t$-element set of $X$ with $|\F_{|T}|>2^t-2^{t-2}$. For any $x,y\in T, x\neq y$, there are $2^{t-2}$ subsets of $T$ that contain $x$, but not $y$. At least one of them must belong to $\F_{|T}$, which proves the claim.  \hfill $\Box$
\medskip

Now we can settle the case $t=3$ in Theorem~\ref{sep}. The first part of Corollary~\ref{cor1} implies that
$$s(n,3)\ge (\lfloor\frac{n}{2}\rfloor+1)(\lceil\frac{n}{2}\rceil+1)+1=\lfloor\frac{n^2}4\rfloor+n+2.$$
On the other hand, the first author~\cite{Fr83} proved that $(n,\lfloor\frac{n^2}4\rfloor+n+2)\rightarrow (3,7)$. Applying Lemma~\ref{lem2} with $m=\lfloor\frac{n^2}4\rfloor+n+2$ and $t=3$, we obtain that $s(n,3)=\lfloor\frac{n^2}4\rfloor+n+2$, which proves part (ii) of Theorem~\ref{sep}.
\medskip

For the rest of the argument, we need some further results from extremal set theory.

Consider again a partition of the $n$-element set $X$ into $t-1$ parts, $X=X_1\sqcup\ldots\sqcup X_{t-1}$. A family $\G\subset\binom{X}{t-1}$ of $(t-1)$-element subsets of $X$ is called a {\em $(t-1)$-uniform hypergraph} or, simply, a {\em $(t-1)$-graph}. If every {\em edge} $E\in \G$ intersects each $X_i$ in precisely $1$ point, then $\G$ is said to be {\em $(t-1)$-partite}.

For graphs, i.e., for $t=3$, Mantel~\cite{Ma07} and Tur\'an~\cite{Tu41} proved that if a graph (2-graph) $\G$ has more than $\lfloor\frac{n^2}{4}\rfloor$ edges (sets), then it contains a {\em triangle}, i.e., there are $x,y,z\in X$ with $\{x,y\}, \{x,z\}, \{y,z\}\in\G$. This bound is best possible, as is shown by a $2$-partite (bipartite) graph whose parts are of size $\lfloor\frac{n}{2}\rfloor$ and $\lceil\frac{n}{2}\rceil$.
\smallskip

For a fixed $t\ge 3$, a {\em generalized triangle} consists of $3$ distinct $(t-1)$-elements sets $E_1, E_2, E_3$ such that $|E_1\cap E_2|=t-2$ and $E_3\supseteq (E_1\setminus E_2)\cup(E_2\setminus E_1)$. For $t=3$, the only generalized triangle is the usual triangle. For $t \ge 4$, however, there are $t-2$ non-isomorphic generalized triangles, depending on the value of $|E_1\cap E_2\cap E_3|$, which can be $0,1,\ldots, t-3$. Obviously, none of these generalized triangles is $(t-1)$-partite.

The Mantel-Tur\'an theorem was extended to $3$-graphs and $4$-graphs by Bollob\'as and Sidorenko, respectively, as follows.

\begin{lemma}\label{bollobas}
Let $X$ be an $n$-element set, $n\ge 4$. For any $k\ge 2$, let $g(n,k)$ denote the maximum size of a $k$-graph $\G$ which does not contain any generalized triangle. Then we have
\smallskip

{\rm (i)}\;\; {\rm (Bollob\'as~\cite{Bo74})}\;\;\;\;\;\;\;
$g(n,3)\le\lfloor\frac{n}{3}\rfloor \lfloor\frac{n+1}{3}\rfloor \lfloor\frac{n+2}{3}\rfloor$;

\smallskip

{\rm (ii)}\; {\rm (Sidorenko~\cite{Si87})}\;\;\;
$g(n,4)\le\lfloor\frac{n}{4}\rfloor \lfloor\frac{n+1}{4}\rfloor \lfloor\frac{n+2}{4}\rfloor \lfloor\frac{n+3}{4}\rfloor$.
\smallskip

\noindent Both results are best possible as is shown by the complete $3$-partite ($4$-partite) $3$-graphs (resp., $4$-graphs) whose parts are as equal as possible.
\end{lemma}

It is not hard to see that part (i) of Lemma~\ref{bollobas} implies the Mantel-Tur\'an theorem for ordinary triangles~\cite{Ka75}. It was proved by Frankl and F\"uredi~\cite{FrF83} (see also \cite{KeM04}) that, if $\G\subset\binom{X}3$ is a 3-graph with $|\G|>g(n,3)$ and $X$ is sufficiently large, then $\G$ also contains a generalized triangle with $E_1\cap E_2\cap E_3=\emptyset$.

One might hope that analogous results hold for $k>4$. However, this is not the case. For $k=5$ and $6$, Frankl and F\"uredi~\cite{FrF89} determined all largest $k$-graphs on $n>n_0$ vertices that do not contain a generalized triangle. These turned out to have  substantially more edges than the balanced complete $k$-partite $k$-graphs, and one can obtain them by ``blowing up'' certain Steiner systems called Witt designs. For a survey on this fascinating problem, consult~\cite{NoY17}.

\smallskip

It remains to establish parts (iii) and (iv) of Theorem~\ref{sep}.
With the notation of Lemma~\ref{bollobas}, they can be rephrased in the following form.
\medskip

\noindent{\bf Theorem 6'} \; {\em Let $n\ge t\ge 4$ and $X=\{1,2,\ldots,n\}$. Let $s=s(n,t)$ denote the smallest number with the property that every family $\F\subset 2^X$ with $|\F|\ge s$ is $t$-separable.  Then we have

$$(\frac{n}{t-1})^{t-1}<s(n,t)\le g(n,t-1)+ 1+\sum_{i=0}^{t-2}\binom{n}{i}.$$

According to Lemma~\ref{bollobas}, for $t=4$ and $5$, the lower bound and the upper bound are asymptotically the same. For $t>5$, the two bounds are asymptotically different, but their order of magnitude is the same, $\Theta(n^{t-1})$.}
\medskip

\noindent{\bf Proof.}  Let $\F\subset 2^X$ be a family satisfying $$|\F|>g(n,t-1)+\sum_{i=0}^{t-2}\binom{n}{i}.$$
We will show that $\F$ is $t$-separable. By Lemma~\ref{lem2}, it is sufficient to prove that there is a $t$-element subset $T\subset X$ such that $|\F_{|T}|\ge 2^{t}-2^{t-2}+1$. According to Lemma~\ref{downward}, we can assume that $\F$ is downward closed.

Therefore, if $\F$ has a member of size at least $t$, then it also has a member $F$ of size precisely $t$. In this case, choosing $T$ to be $F$, we have $|\F_{|T}|=2^t$, and we are done.

Thus, we can assume that the number of $(t-1)$-element members in $\F$ is larger than $g(n,t-1)$.
By the definition of $g(n,t-1)$, the $(t-1)$-graph $\G\subset\F$ formed by these edges contains a generalized triangle $E_1, E_2, E_3$. We can assume without loss of generality that $E_1=\{1,2,\ldots,t-2, t-1\}, E_2=\{1,2,\ldots,t-2,t\},$ and $\{t-1,t\}\subset E_3$. Set $T=\{1,2,\ldots, t\}$. Then we have $|T|=t$ and
$\G_{|T}\supset (2^{E_1}\cup 2^{E_2})\sqcup\{t-1,t\}.$ As
$$|2^{E_1}\cup 2^{E_2}|=2^{|E_1|}+2^{|E_2|}-2^{|E_1\cap E_2|}=2^{t-1}+2^{t-1}-2^{t-2}=2^t-2^{t-2},$$
we obtain that $|\F_{|T}|\ge|\G_{|T}|\ge 2^t-2^{t-2} +1$, as required. This completes the proof of the upper bound in Theorem~6'. The lower bound is given by  Corollary~\ref{cor1}. \hfill $\Box$

\section{Open problems, concluding remarks}

It would be interesting to close the gaps between the lower and upper bounds in Theorems~\ref{thm:main} and~\ref{sep}.
\smallskip

{\bf 4.1.}  What happens if, instead of concentrating on maximal intersecting families of $\frac{n}2$-element subsets of an $n$ element set, as we did in Theorem~\ref{thm:main}, we consider {\em all} maximal intersecting families of subsets of $X$, with no restriction on the sizes of the subsets? In particular, we can ask the following.

\begin{problem}\label{problem1}
Determine or estimate the largest integer $k^*=k^*(n)$ such that for every maximal intersecting family $\F$ of subsets of an $n$-element set, one can find a shattered matching of size $k^*$.
\end{problem}

{\bf 4.2.} Given a family $\F\subset 2^X$ and a system $\M$ of pairwise disjoint {\em $r$-element} subsets of $X$ for some $r\ge 3$, we say that $\M$ is {\em shattered} by $\F$ if no matter how we pick one element from each $r$-tuple of $\M$, there is a member $F\in\F$ which carves out precisely these elements of $\cup\M$.

\begin{problem}\label{problem0}
For any even integer $n\ge 4$ and $r\ge 3$, determine or estimate $k_r(n)$, the largest integer $k$ such that for every maximal intersecting family $\F$ of $\frac{n}2$-element subsets of an $n$-element set, one can find a system $\M$ of $k$ pairwise disjoint $r$-element sets which is shattered by $\F$.
\end{problem}

We trivially have $k_r(n)<\frac{n}{2(r-1)}$, whenever $n$ is a multiple of $2(r-1)$. If we take $\frac{n}{2(r-1)}$  $r$-tuples and pick one element from each, then there is a unique way how to add further elements from the remainder to obtain an $\frac{n}2$-element set.
\smallskip

As in Problem~\ref{problem1}, here we can also relax the condition that every member of our maximal intersecting family is of size $\frac{n}2$. Furthermore, in the spirit of Problem~\ref{problem2}, we can completely drop the restriction that $\F$ is intersecting, and we can ask {\em how large} $\F$ needs to be in order to ensure that it shatters some system of $k$ pairwise disjoint $r$-element sets.
\medskip

{\bf 4.3.}   Following \cite{FrP84}, we call a family of $t$ sets $F_1,\ldots,F_t\subset X$  {\em disjointly representable} if there exist $x_1,\ldots,x_t\in X$ with the property that $x_i\in F_j$ if and only if $i=j$. In other words, a family of sets is disjointly representable if and only if none of its members is completely covered by the union of the others.

Modifying the question addressed in Theorem~\ref{sep}, we can ask the following.

\begin{problem}\label{disrep}
Let $n\ge t\ge 2$. Determine or estimate the smallest number $r=r(n,t)$ with the property that every family $\F$ of at least $r$ subsets of an $n$-element set has $t$ disjointly representable members.
\end{problem}

It follows from the definition that if a family has $t$ disjointly representable members, then it is $t$-separable. Therefore, we have $r(n,t)\ge s(n,t)$ for every $n$ and $t$.

The proof of part (i) of Theorem~\ref{sep} also gives $r(n,2)=s(n,2)=n+2$.

\begin{claim}
$r(n,3)=\binom{n}2+n+2>s(n,3)$ for every $n\ge 3$.
\end{claim}

\noindent{\bf Proof.} The upper bound immediately follows from Theorem~\ref{sauer}.

To prove the lower bound, consider the following subsets of $X=\{1,2,\ldots,n\}$. For any two elements $a<b$ of the auxiliary set $A=\{\frac12, 1+\frac12,\ldots,n+\frac12\}$, let
$$F(a,b)=\{x\in X\; :\; x<a\}\cup\{x\in X\;|\;x>b\}.$$
The family $\F=\{F(a,b) : a,b\in A \mbox{ and } a<b\}\cup\{\emptyset\}$ has $\binom{n+1}2+1=\binom{n}2+n+1$ members, and it has no 3 disjointly representable members. Indeed, given any 3 elements $x_1<x_2<x_3\in X$, there is no $F\in\F$ with $F\cap\{x_1,x_2,x_3\}=\{x_2\}$.  \hfill $\Box$

\end{document}